\documentclass[10pt,reqno]{amsart}
\usepackage{amssymb}
\usepackage{amsfonts}
\usepackage[all]{xy}
\usepackage[english]{babel}

\def\barM {{\overline M}}

\begin{document}

\title{Scrambled and distributionally scrambled $n$-tuples}

\author{Jana Dole\v zelov\'a*}
 \thanks{*Mathematical Institute, Silesian University, CZ-746 01 
Opava, Czech Republic\\ \emph{E-mail address:} jana.dolezelova@math.slu.cz}

\address{J . Dole\v zelov\'a, Mathematical Institute, Silesian University, CZ-746 01 
Opava, Czech Republic}

\email{jana.dolezelova@math.slu.cz, }

\thanks{The research was supported by grant SGS/2/2013 from  the Silesian University in Opava. Support of this institution is gratefully acknowledged.}

\begin{abstract} 
This article investigates the relation between the distributional chaos and the existence of a scrambled triple. We show that for a continuous mapping $f$ acting on a compact metric space $(X,d)$, the possession of an infinite extremal distributionally scrambled set is not sufficient for the existence of a scrambled triple. We also construct an invariant Mycielski set  with an uncountable extremal distributionally scrambled set without any scrambled triple.\\
{\small {2000 {\it Mathematics Subject Classification.}}
Primary 37D45; 37B10.
}
\end{abstract}

\maketitle
\pagestyle{myheadings}
\markboth{Jana Dole\v zelov\' a}
{Scrambled and distributionally scrambled $n$-tuples}
\section{Introduction}
The first  definition of chaotic pairs appeared in the paper \cite{LY} by Li and Yorke in 1975. One of the most important extensions of the concept of Li-Yorke chaos is distributional chaos introduced in \cite{SchSm}. This extended definition is much stronger - there are many mappings which are chaotic in the sense of Li-Yorke but not distributionally chaotic. Another way how to extend the Li-Yorke chaos is looking on dynamics of tuples instead dynamics of pairs. Since Xiong \cite{X} and Sm\'ital \cite{S} constructed some interval maps with zero topological entropy which are Li-Yorke chaotic, the Li-Yorke chaos is not sufficient condition for positive topological entropy. But interval maps with zero topological entropy never contain scramled triples \cite{L} and hence existence of scrambled triple implies positive topological entropy. Consequently we can find a dynamical system which is Li-Yorke chaotic but contains no scrambled triple. The natural question was if there is a dynamical system which is distributionally chaotic but contains no distributionally scrambled triple. Example in \cite{O} contains no distributionally scrambled triple but still there were some scrambled triples (in the sense of Li-Yorke) and therefore another open problem appeared - is there a distributionally chaotic dynamical system without any scrambled triple? In this paper, we construct a dynamical system which possess an infinite extremal distributionally scrambled set but without any scrambled triple. We show the existence of an invariant Mycielski set which possess an uncountable extremal distributionally scrambled set and has no scrambled triple. Because some scrambled triples occur in the closure of this invariant set, the following question remains open:
\emph{Does the existence of uncountable distributionally scrambled set imply the existence of a scrambled triple?}
\section{Terminology}
Let $(X,d)$ be a non-empty compact metric space.  Let us denote by $(X,f)$ the \emph{topological dynamical system}, where $f$ is a continous self-map acting on $X$.  We define the \emph{forward orbit} of $x$, denoted by $Orb^+_f(x)$ as the set $\{f^n(x):n\geq0\}$. A non-empty closed invariant subset
$Y\subset X$ defines naturally a subsystem $(Y, f )$ of $(X, f )$. For $n\geq 2$, we denote by $(X^n,f ^{(n)})$ the product system $(X\times X\times...\times X,f\times f\times...\times f)$ and put $\Delta^{(n)} = \{(x_1, x_2, . . . , x_n)\in X^n : x_i =
x_j \mbox{ for some } i\neq j\}$. By a \emph{perfect set} we mean a nonempty compact set without isolated points. A \emph{Cantor set} is a nonempty, perfect and totally disconnected set.  A set $D\subset X$ is \emph{invariant} if $f(D)\subset D$. A \emph{Mycielski set} is defined as a countable union of Cantor sets.\\

\paragraph{\bf Definition 1} A tuple $(x_1, x_2,...,x_n)\in X^n$ is called \emph{$n$-scrambled} if 
\begin{equation}\liminf_{k\to\infty} \max_{1\leq i<j\leq n} d(f^k(x_i),f^{k}(x_j))=0\end{equation} and \begin{equation}\limsup_{k\to\infty} \min_{1\leq i<j\leq n} d(f^k(x_i),f^{k}(x_j))>0.\end{equation}
A subset $S$ of $X$ is called \emph{$n$-scrambled} if every $n$-tuple $(x_1,x_2,...,x_n)\in S^n\setminus \Delta^{(n)}$ is $n$-scrambled. The system $(X,f)$ is called \emph{$n$-chaotic} if there exists an uncountable $n$-scrambled set.\\

\paragraph{\bf Definition 2} For an $n$-tuple $(x_1, x_2,...,x_n)$ of
points in $X$, define the \emph{lower distribution function} generated by $f$ as
$$\Phi_{(x_1, x_2,...,x_n)}(\delta)=\displaystyle\liminf_{m\to\infty}\frac{1}{m}\#\{0<k<m;\min_{1\leq i<j\leq n} d(f^k(x_i),f^{k}(x_j))<\delta\},$$
and the \emph{upper distributional function} as 
$$\Phi^*_{(x_1, x_2,...,x_n)}(\delta)=\displaystyle\limsup_{m\to\infty}\frac{1}{m}\#\{0<k<m;\max_{1\leq i<j\leq n} d(f^k(x_i),f^{k}(x_j))<\delta\},$$
where $\#A$ denotes the cardinality of the set $A$.\\ 
A tuple $(x_1, x_2,...,x_n)\in X^n$ is called \emph{distributionally $n$-scrambled} if 
$$\Phi^*_{(x_1, x_2,...,x_n)}\equiv 1 \mbox{  and  } \Phi_{(x_1, x_2,...,x_n)}(\delta)=0, \mbox{  for some  }0<\delta\le \text{diam }X .
$$
A subset $S$ of $X$ is called \emph{distributionally $n$-scrambled} if every $n$-tuple $(x_1,x_2,...,x_n)\in S^n\setminus \Delta^{(n)}$ is distributionally $n$-scrambled. The
system $(X,f)$ is called \emph{distributionally $n$-chaotic} if there exists an uncountable distributionally $n$-scrambled set.\\

\paragraph{\bf Definition 3} A subset $S$ of $X$ is called \emph{extremal distributionally $n$-scrambled} if every $n$-tuple $(x_1,x_2,...,x_n)\in S^n\setminus \Delta^{(n)}$ is distributionally $n$-scrambled with $\Phi_{(x_1, x_2,...,x_n)}(\delta)=0$, for any $\delta<\text{diam }X.$\\

Let $A=\{0,1,...,n-1\}, n\geq 2$, be a finite alphabet and $\Sigma_n$ a set of all infinite sequences on A, that is, for
$u\in \Sigma_n$, $u=u_1 u_2 u_3\ldots$, where $u_i\in A$ for all $i\geq 1$.
We define a metric on $\Sigma_n$ by
$$d(u,v)=\sum_{i=1}^{\infty} \frac{\delta(u_i,v_i)}{2^i},$$
where 
 $$\delta(u_i,v_i) = \left\{ \begin{array}{r@{\quad}c}
    0, &  u_i=v_i\\
    1, & u_i\neq v_i.\\ \end{array} \right.
    $$ 
The \emph{shift} transformation is a continuous map $\sigma:\Sigma_n\rightarrow \Sigma_n$ given by $\sigma(u)_i=u_{i+1}$.The dynamical system $(\Sigma_n,\sigma)$ is called the one-sided shift on $n$ symbols. 
Any closed subset $X\subset\Sigma_n$ invariant for $\sigma$ is called a subshift of $(\Sigma_n,\sigma)$.\\
Any finite string $B$ of some $u\in\Sigma_n$ is called a word (or a block) and the length of $B$ is denoted by $|B|$. Let $B=b_1\ldots b_n$ and $G=g_1\ldots g_n$ be words. Denote by $BG=b_1\ldots b_n g_1\ldots g_n$ and for the case of $\Sigma_2$ denote by  $\bar{B}$ the binary complement of $B$.\\ 

The \emph{Morse block} $M_i$ is defined inductively such that $M_0=0$, and $M_i=M_{i-1}\barM_{i-1}$, for all $i>0$. The \emph{Morse sequence} $m\in\Sigma_2$ is the limit of the Morse blocks, i.e. $m=\lim_{i\to\infty} M_i$ . This sequence $m$ generates the infinite Morse minimal set $M=cl\{ m,\sigma(m),\sigma^2(m),\ldots \}$ and it is known that, for all words $B\subset m$, the sequence $m$ contains no block $BBb$, where $b$ is the first element of block $B$ (cf. \cite{GH}). Denote this property by $\mathcal{P}$.

\section{Scrambled and distributionally scrambled $n$-tuples}

We will show that the existence of an infinite extremal distributionally  scrambled set doesn't imply the existence of a scrambled triple. Then we construct an invariant Mycielski set with an uncountable extremal distributionally scrambled
set without any scrambled triple. \\
\paragraph{\bf Lemma 1}
\emph{Let $\{a_i\}^{\infty}_{i=1}$ be a strictly increasing sequence of positive integers such that, for every $n\geq 1$, $a_n$ and $n$ have the same parity. Then the point $x=M_{a_1}M_{a_2}M_{a_3}\ldots $ is contained in the Morse minimal set.}
\begin{proof}
Since $a_n$ and $a_{n+1}$ have different parity, for every $n\geq 1$, the Morse block $M_{a_{n+1}}$ ends with $\bar M_{a_{n}}$ 
$$M_{a_{n+1}}=M_{a_{n}}\bar M_{a_{n}}\bar M_{a_{n}}M_{a_{n}}\ldots \bar M_{a_{n}}.$$
By the construction of the Morse sequence $m$ and the previous property, we can observe
$$\sigma^{3\cdot 2^{a_2}-2^{a_1}}(m)=M_{a_1}M_{a_2}\bar M_{a_2}M_{a_2}M_{a_2}\bar M_{a_2}\ldots$$
$$\sigma^{3\cdot 2^{a_3}-2^{a_1}-2^{a_2}}(m)=M_{a_1}M_{a_2}M_{a_3}\bar M_{a_3}M_{a_3}M_{a_3}\bar M_{a_3}\ldots$$
$$\vdots$$
Since $\sigma^{r_n}(m)$ starts with blocks $M_{a_1}M_{a_2}M_{a_3}\ldots M_{a_n}$, where $r_n=3\cdot2^{a_n}-\sum_{i=1}^{n-1}2^{a_i}$, for all $n>1$, the point $x=\lim_{n\to\infty}\sigma^{r_n}(m)$ is contained in the Morse minimal set.
\end{proof}
\paragraph{\bf Theorem 1}
\emph{There exists a dynamical system $X$ with an infinite extremal distributionally scrambled set but without any scrambled triple.}

\begin{proof}
Let $W_1,W_2,W_3,\ldots$ be the following infinite decomposition of even numbers into infinite sets:
$$W_1=\{ 2^n\cdot1,n\geq 1\},$$
$$W_2=\{ 2^n\cdot3,n\geq 1\},$$
$$W_3=\{ 2^n\cdot5,n\geq 1\},$$
$$\vdots$$
Let $\{a_n\}_{n=1}^{\infty}$ be an increasing sequence of positive integers with $\lim_{n\to\infty}a_n/a_{n+1}=0$ and, for every $n\geq 1$, $a_n$ and $n$ have the same parity. Then \begin{equation}\lim_{n\to\infty}\frac{\sum_{i=1}^{n-1} 2^{a_i}}{2^{a_n}}=0.\end{equation} 
We construct the point $x^i$ as a sequence of blocks $$M_{a_1}M_{a_2}^iM_{a_3}M_{a_4}^iM_{a_5}M_{a_6}^iM_{a_7}M_{a_8}^i\ldots$$
where $$M_{a_j}^i =
\begin{cases}
 M_{a_j}, & \text{if }j\notin W_i\\
\bar{M}_{a_j}, & \text{if }j\in W_i.
\end{cases}$$
{\bf Remark 1} Let $i$ be a fixed positive integer. Then the first complementary block $\bar M_{a_j}$ appears in the construction of $x^i$ for $j=2\cdot (2i-1),$
$$x^i=M_{a_1}M_{a_2}M_{a_3}M_{a_4}M_{a_5}\ldots M_{a_{2\cdot (2i-1)-1}}\bar M_{a_{2\cdot (2i-1)}}\ldots.$$
 Hence the sequence $\{x^i\}_{i=1}^{\infty}$ converges and $\lim_{i\to\infty}x^i=x$, where $x$ is the point constructed in Lemma 1.\\
Let $D=\{x^i\}_{i=1}^{\infty}.$ We claim $D$ is distributionally 2-scrambled set and $X=cl(\cup_{i=0}^{\infty}\sigma^i(D))$ is the wanted dynamical system.\\

\emph{I. $D$ is extremal distributionally 2-scrambled set}\\
Let $(x^i,x^j)\in D^2$ be a pair of distinct points. For simplicity denote $s_k=\sum_{n=1}^k 2^{a_n}$. Let $l$ be a fixed positive integer and $\epsilon=\frac{1}{2^l}$. Since $(x^i)_n=(x^j)_n$ if $s_{2k}<n\leq s_{2k+1}$, for any $k>0$, we have $d(\sigma^n(x^i),\sigma^n(x^j))<\epsilon$ for all $s_{2k}<n<s_{2k+1}-l$. By (3), $$\displaystyle\lim_{k\to\infty}\frac{2^{a_{2k+1}}-l}{2^{a_{2k+1}}+s_{2k}}=1,$$  so it is easy to see that $\Phi^*_{x^i,x^j}(\epsilon)=1$, for arbitrary small $\epsilon$, and hence $\Phi^*_{x^i,x^j}\equiv 1$.\\
On the other hand, there is a sequence $\{l_k\}_{k=1}^{\infty}\subset W_i\cup W_j$ such that $(x^i)_n=(\overline{{x}^j})_n$  if $s_{l_k-1}<n\leq s_{l_k}$, for any integer $k$. Since $(x^i)_m=(\overline{{x}^j})_m$, for $m=1,2,\ldots r$, implies $d(x^i,x^j)\geq \sum_{m=1}^r\frac{1}{2^m}$ and $(\sigma^n(x^i))_m=(\overline{\sigma^n({x}^j)})_m$, for  all $s_{l_k-1}<n<s_{l_k}-r$ and $m=1,2,\ldots r$, it follows $d(\sigma^n(x^i),\sigma^n(x^j))\geq \sum_{m=1}^r\frac{1}{2^m}$ for all $s_{l_k-1}<n<s_{l_k}-r$. Because $$\lim_{k\to\infty}\frac{s_{l_k-1}}{{s_{l_k-1}+2^{a_{l_k}}-r}}=0,$$  it is easy to see that $\Phi_{x^i,x^j}(\sum_{m=1}^r\frac{1}{2^m})=0$, for arbitrary large $r$, and hence $\Phi_{x^i,x^j}(\delta)=0$, for any $0<\delta<1.$\\

\emph{II. $\bigcup_{i=0}^{\infty}\sigma^i(D)$ has no scrambled triples}\\
Let $(x^i,x^j,x^k)\in D^3\setminus \Delta^{(3)}$. Since $$M_{2n}^{i}=M_{2n}^{j} \mbox{ or } M_{2n}^{i}=M_{2n}^{k }\mbox{ or } M_{2n}^{j}=M_{2n}^{k},$$ and
$M_{2n-1}$ is the common block for all $x^i,x^j,x^k$ and every $n>1$, we can assume $$\lim_{n\to\infty} \min \{ d(\sigma^n(x^{i}),\sigma^{n}(x^{j})), d(\sigma^n(x^{i}),\sigma^{n}(x^{k})), d(\sigma^n(x^{k}),\sigma^{n}(x^{j}))\}=0$$ and consequently, condition (2) is not satisfied and $D$ has no scrambled triples. For the same reason  $\sigma^n(D)$ has no scrambled triples, for any $n>0$. It follows that any potential scrambled triple in $X$ must contain some pair $\sigma^p(u),\sigma^q(v)$, where $p<q$ and $u,v\in D$. To prove that for such tuple the condition (1) is not fulfilled, it  is sufficient to show that
$$\liminf_{k\to \infty} d(\sigma^k(\sigma^p(u)),\sigma^k(\sigma^q(v)))>0,$$ where $p<q$ and $u,v\in D.$
Assume the contrary - let $\liminf_{k\to \infty} d(\sigma^k(\sigma^p(u)),\sigma^k(\sigma^q(v)))=0$ and denote $r=q-p>0$. Then we can find an infinite subsequence $\{k_n\}_{n=1}^{\infty}$ such that both $\sigma^{k_n}(\sigma^p(u))$ and $\sigma^{k_n}(\sigma^q(v))$ begin with the same block $G_n$ of length $14r$ and obviously these blocks can be found also in the sequence $u$ and, shifted by $r$, in $v$. For sufficiently large $n$, $G_n$ is in the sequence $u$ either contained in some Morse block or is on the edge of two following Morse blocks, but at least the first $7r$ digits or the last $7r$ digits of block $G_n$ are contained in a single Morse block. Denote this block $M^{(u)}_{a_j}$, where $M^{(u)}_{a_j}$ is either $M_{a_j}$ or $\overline{M}_{a_j})$, depending on $u$ and $j$, and these $7r$ consecutive digits of $G_n$ by $G=g^u_1g^u_2...g^u_{7r}$. This $G$ appears in $v$ shifted by $r$, so we can conclude that the first $6r$ digits $g^v_1g^v_2...g^v_{6r}$ of $G$ in $v$ are in $M^{(v)}_{a_j}$. The block $M^{(v)}_{a_j}$ is either the same Morse block as $M^{(u)}_{a_j}$ or its binary complement. In the first case, $g_1^u=g^v_1=g^u_{r+1}=g^v_{r+1}=g^u_{2r+1}=...=g^u_{5r+1}$ and similarly for $g_2^u,...,g_r^u$ and therefore we obtained a block BBBBBB which is impossible since $M^{(u)}_{a_j}$ is a Morse block. In the second case, $g_1^u=g^v_1=\overline{g}^u_{r+1}=\overline{g}^v_{r+1}=\overline{\overline{g}}^u_{2r+1}=g^u_{2r+1}=...=g^u_{4r+1}$ and similarly for $g_2^u,...,g_{2r}^u$ and therefore we obtained a block BBB which is a contradiction with $M^{(u)}_{a_j}$ is a Morse block.\\
{\bf Remark 2} By the same argument, $(x^i,x^j,x)$ is never scrambled triple, where $x=\lim_{i\to\infty}x^i$.

\emph{III. If $y\in X\setminus \bigcup_{i=0}^{\infty}\sigma^i(D)$, then there exists a nonnegative integer $n$ such that $\sigma^n(y)$ is contained in the Morse minimal set $M$.}\\
Suppose the contrary. By \cite{GH}, the points of Morse minimal sets are characterised by the property  $\mathcal{P}$ and therefore we can find two distinct blocks $B_1B_1b_1$ and $B_2B_2b_2$ which appears in $y$, where $b_1$ denotes the first element of $B_1$ and $b_2$ denotes the first element of $B_2$. Suppose that  the last element of $B_1B_1b_1$ is on the $k_1$-th position in $y$, the last element of $B_2B_2b_2$ is on the $k_2$-th position in $y$ and $k_2>k_1+2\cdot|B_2|$. Since $y$ is contained in the closure of  $\cup_{i=0}^{\infty}\sigma^i(D)$, there are sequences $\{x^{m_{i}}\}_{i=1}^{\infty}$ and $\{n_i\}_{i=1}^{\infty}$ such that
$$y=\lim_{i\to\infty}\sigma^{n_i}(x^{m_i}),$$ and suppose all sequences $\sigma^{n_i}(x^{m_i})$ have the same first $k_2$ symbols. Let $j$ be an integer such that $2^{a_j}>k_2$. We can observe that $n_i$ is bounded by $\sum_{l=1}^{j-1}2^l$ for all $i>0$ - for larger $n_i$ the block of the first $k_2$ symbols would be part of a single Morse block, but by assumption there is $B_1B_1b_1$ inside of the block of first $k_2$ symbols. Hence there is a nonnegative integer $N$ and a subsequence $\{x^{m_{i_k}}\}_{k=1}^{\infty}$ such that $y=\lim_{i\to\infty}\sigma^{n_i}(x^{m_i})=\lim_{i\to\infty}\sigma^N(x^{m_{i_k}})=\sigma^N(x)$, where the last identity follows from the Remark 1. By Lemma 1, $y=\sigma^N(x)\in M$ and this is a contradiction with assumptions. 

\emph{IV. $X$ has no scrambled triples}
By \cite{G} the Morse minimal set is a distal system, i.e. all pairs in this set are either distal or asymptotic. Hence the only potential scrambled triples are $(x^i,x^j,y)$, where $y\in X\setminus \cup_{i=0}^{\infty}\sigma^i(D)$ and $x^i,x^j\in D$. By the previous step, it is sufficient to show that $(x^i,x^j,y)$, where $y\in M$, is not a scrambled triple. Since $y\in M$ and $x\in M$ by Lemma 1, the pair $(x,y)$ is either distal or asymptotic:\\
a) $(x,y)$ is distal pair\\
Sequences $x,x^i,x^j$ are exactly the same except of blocks $M^i_{a_l}$ and $M^j_{a_l}$ where $l\in W_i\cup W_j$ and it holds $M^i_{a_l}=\overline{M}^j_{a_l}$, for $l\in W_i\cup W_j$. Therefore $(x^i,x^j,y)$ is not proximal.\\
b) $(x,y)$ is asymptotic pair\\
The triple $(x^i,x^j,x)$ is not scrambled by Remark 2, therefore $(x^i,x^j,y)$ is not scrambled.

\end{proof}
To prove Theorem 2, we need the next lemma:\\

\paragraph{\bf Lemma 2}
\emph{There is a Cantor set $B\subset \{0,1\}^{\mathbb{N}}$ such that, for any distinct $\alpha=\{\alpha(i)\}^{\infty}_{i=1}$ and $\beta=\{\beta(i)\}^{\infty}_{i=1}$ in $B$, the set 
\begin{equation}\label{eq:lm1}
\{j\in\mathbb{N};\alpha(j)\neq\beta(j)\}\hspace{.2cm}\mbox{is infinite.}
\end{equation}}

\begin{proof}
By Lemma 5.4 in \cite{SchSm} there is an uncountable Borel set $B\subset \{0,1\}^{\mathbb{N}}$ satisfying (\ref{eq:lm1}). The result follows from Alexandrov-Hausdorff Theorem \cite{K}.
\end{proof}

\paragraph{\bf Theorem 2}
\emph{There exists an invariant Mycielski set $X\subset\Sigma_2$ with an uncountable extremal distributionally
2-scrambled set but without any 3-scrambled tuple.}\\

\begin{proof}

We will denote by $M^0_i$ the Morse block $M_i$ and by $M^1_i$ the binary complement of $M_i$.  Let $\alpha=\{\alpha(i)\}^{\infty}_{i=1}$ be a point of $B$ where $B$ is the set from Lemma 2. Let $\{a_n\}_{n=1}^{\infty}$ be an increasing sequence of positive integers with $\lim_{n\to\infty}a_n/a_{n+1}=0$. Then \begin{equation}\lim_{n\to\infty}\frac{\sum_{i=1}^{n-1} 2^{a_i}}{2^{a_n}}=0.\end{equation} We construct a point $x^{\alpha}$ as a sequence of blocks 
$$M^{\alpha_1}_{a_1}M^{0}_{a_2}M^{\alpha_2}_{a_3}M^{0}_{a_4}M^{\alpha_3}_{a_5}M^{0}_{a_6}M^{\alpha_4}_{a_7}\ldots .$$
Let $D=\{x^{\alpha};\alpha\in B\}$. We claim $D$ is extremal distributionally 2-scrambled set and $X=\cup_{i=0}^{\infty} \sigma^i(D)$ is the wanted Mycielski set. 

\emph{ I. $D$ is extremal distributionally 2-scrambled set}\\
Let $(u,v)\in D^2$ be a pair of distinct points. For simplicity denote $s_i=\sum_{j=1}^i 2^{a_j}$. Let $l$ be a fixed integer and $\epsilon=\frac{1}{2^l}$. Since $u_i=v_i$ if $s_{2k-1}<i\leq s_{2k}$, for any $k>0$, we have $d(\sigma^i(u),\sigma^i(v))<\epsilon$ for all $s_{2k-1}<i<s_{2k}-l$. By (4), $$\displaystyle\lim_{k\to\infty}\frac{2^{a_{2k}}-l}{2^{a_{2k}}+s_{2k-1}}=1,$$  so it is easy to see that $\Phi^*_{u,v}(\epsilon)=1$, for arbitrary small $\epsilon$, and hence $\Phi^*_{u,v}\equiv 1$.\\
On the other hand, there is a sequence $\{l_k\}_{k=1}^{\infty}$ such that $u_i=\bar{v}_i$  if $s_{l_k-1}<i\leq s_{l_k}$, for any integer $k$. Since $u_m=\bar{v}_m$, for $m=1,2,\ldots r$, implies $d(u,v)\geq \sum_{m=1}^r\frac{1}{2^r}$ and it follows $d(\sigma^i(u),\sigma^i(v))\geq \sum_{m=1}^r\frac{1}{2^r}$ for all $s_{l_k-1}<i<s_{l_k}-r$ and $m=1,2,\ldots r$. Because $$\lim_{k\to\infty}\frac{s_{l_k-1}}{{s_{l_k-1}+2^{a_{l_k}}-r}}=0,$$  it is easy to see $\Phi_{u,v}(\sum_{m=1}^r\frac{1}{2^m})=0$, for arbitrary large $r$, and hence $\Phi_{u,v}(\delta)=0$, for any $0<\delta<1.$\\

\emph{ II. $X$ has no scrambled triples}\\
Let $(x^\alpha,x^\beta,x^\gamma)\in D^3\setminus \Delta^{(3)}$. Since $\alpha_i,\beta_i,\gamma_i\in\{0,1\}$, for any integer $i$,  $$M_{2i-1}^{\alpha_i}=M_{2i-1}^{\beta_i} \mbox{ or } M_{2i-1}^{\alpha_i}=M_{2i-1}^{\gamma_i }\mbox{ or } M_{2i-1}^{\beta_i}=M_{2i-1}^{\gamma_i },$$ and
$M_{2i}^{0}$ is the common block for all $x^\alpha,x^\beta,x^\gamma$, we can assume $$\lim_{k\to\infty} \min \{ d(\sigma^k(x^{\alpha}),\sigma^{k}(x^{\beta})), d(\sigma^k(x^{\alpha}),\sigma^{k}(x^{\gamma})), d(\sigma^k(x^{\gamma}),\sigma^{k}(x^{\beta}))\}=0$$ and consequently, condition (2) is not satisfied and $D$ has no scrambled triples. For the same reason  $\sigma^i(D)$ has no scrambled triples, for any $i>0$. It follows that any potential scrambled triple in $X$ must contain some pair $\sigma^p(u),\sigma^q(v)$, where $p<q$ and $u,v\in D$. The fact that for such tuple the condition (1) is not fulfilled, can be proven in the similar way to the second step of the proof of the previous theorem.\\

\emph{III. $X$ is a Mycielski set}\\
 Let $h:B\rightarrow D$ be a bijection such that, for all $\alpha\in B$,$$h(\alpha)=x^{(\alpha)}.$$ 
To prove that $h$ is homeomorphism, it is sufficient to show that $h$ is continuous. Let $\{ \alpha_m\}_{m=1}^{\infty}$ be a converging sequence in $B$, i.e. $\lim_{m\rightarrow\infty} \alpha_m =\alpha$. Then for an arbitrary $i>0$ there is an $m_0$ such that, for all $m>m_0$, the first $i$ members of the sequences $\alpha_m$ and $\alpha$ are equal. Therefore also the first $2^{a_1}+2^{a_1}+2^{a_2}+\ldots +2^{a_{(2i-1)}}$ members of $x^{\alpha_m}$ and $x^{\alpha}$ are equal and this exactly means $\lim_{m\rightarrow\infty}x^{ \alpha_m} =x^{\alpha}$, hence $h$ is homeomorphism and $D$ is a Cantor set. Since $D$ is distributionally 2-scrambled, the mapping $\sigma^i|_D:D\rightarrow \sigma^i(D)$ is one-to-one and $\sigma^i|_D$ is homeomorphism for every $i\geq 1$. Thus $X$ is the union of Cantor sets.
\end{proof}
{\bf Remark 3} There are some scrambled triples $(x^{\alpha},x^{\beta},x)$ in the closure of $X$, where $x\in cl(X)\setminus X$ and depends on the parity of $\{a_n\}_{n=1}^{\infty}$. If $a_n$ and $n$ have the same parity, then there exist $\alpha$ and $\beta$ such that $(x^{\alpha},x^{\beta},x)$ is scrambled, where $x=\overline{M}_{a_1}\overline{M}_{a_2}\overline{M}_{a_3}\overline{M}_{a_4}\ldots$.\\

\paragraph{\bf Acknowledgment}
I sincerely thank my supervisor, Professor Jaroslav Sm\' ital, for valuable guidance. I~am grateful for his constant support and help.\\


\begin{thebibliography}{m}

\bibitem[1]{SchSm}{\sc Schweizer B., Sm\'ital J.}, {\em Measures of chaos and a spectral decomposition of dynamical systems on the interval\/}, Trans. Amer. Math. Soc. {\bf 344}, (1994), 737 -- 754.         

\bibitem[2]{K}{\sc Kuratowski K.}, {\em Topologie}, Vol. II., Academic Press and Polisch Scientific Publischers, (1968).

\bibitem[3]{GH}{\sc Gottschalk W. H., Hedlund G. A.}, {\em A characterization of the Morse minimal set\/}, Proc. Amer. Math. Soc. {\bf 15}, (1964), 70--74.
\bibitem[4]{LY}{\sc Li T., Yorke J.}, {\em Period three implies chaos\/}, Amer. Math. Monthly {\bf 82}, (1975), 985--992.         
\bibitem[5]{X}{\sc Xiong J.}, {\em A chaotic map with topological entropy zero\/}, Acta Math. Sci. {\bf 6}, (1986), 439--443.   
\bibitem[6]{L}{\sc Li J.}, {\em Chaos and entropy for interval maps\/}, J. Dynam. Differential Eq. {\bf 23}, (2011), 333--352.   
\bibitem[7]{S}{\sc Sm\'ital J.}, {\em Chaotic functions with zero topological entropy\/}, Trans. Amer. Math. Soc. {\bf 297}, (1986), 269 -- 282.   
\bibitem[8]{O}{\sc Oprocha P., Li J.}, {\em On $n$-scrambled tuples and distributional chaos in a sequence\/}, J. Difference Eq. Appl. {\bf 19}, (2013), 927--941
 \bibitem[9]{G}{\sc Blanchard F., Glasner E., Kolyada S., Maas A.}, {\em On Li-Yorke pairs}, J. reine angew. Math. {\bf 547}, (2002), 51--68
     \end{thebibliography}
\end{document}